\documentclass[11pt,oneside]{amsart}

\usepackage{mathtools,amsthm,amssymb,amsfonts,amsmath}
\usepackage{graphicx,txfonts}
\usepackage{hyperref}
\usepackage{verbatim}
\usepackage{csquotes}
\usepackage{geometry}
\usepackage{color}
\usepackage{tikz-cd}
\usetikzlibrary{cd}
\usepackage{setspace}
\usepackage{booktabs}
\usepackage{todonotes}
\usepackage[all]{xy}
\usepackage[shortlabels]{enumitem}
\usepackage[capitalise]{cleveref}
\usepackage{textgreek}
\usepackage{dutchcal}
\usepackage{pifont}
\usepackage{charter,eulervm}
\usepackage[mathscr]{euscript}

\DeclareMathOperator{\coz}{Coz}
\DeclareMathOperator{\bcoz}{bCoz}

\newcommand{\map}[3]{ #1 \colon #2 \to #3 }

\newcommand\DDarrow{\mathrel{\rotatebox[origin=c]{-90}{$\twoheadrightarrow$}}}

\def\O{\mathcal{O}}

\newcommand{\N}{\mathbb{N}}

\newcommand{\Rn}{\mathbb{R}}

\newcommand{\Coz}{\mathord{\textsf{Coz}}}

\newcommand{\rz}{\mathord{\textsf{rZ}}}

\newcommand{\ignore}[1]{}

\newcommand{\cb}{\prec\!\!\prec}

\newcommand{\meet}{\wedge}
\newcommand{\join}{\vee}

\newcommand{\Ss}{\mathsf{S}}

\newcommand{\tbigvee}{\mathop{\textstyle \bigvee}}
\def\open{{\mathfrak o}} 
\def\close{{\mathfrak c}}

\newcommand{\zs}{\mathsf{ZS}}
\newcommand{\cs}{\mathfrak{c}}
\newcommand{\rzs}{\mathsf{rZS}}
\newcommand{\bzs}{\mathsf{bZS}}
\newcommand{\cozs}{\mathsf{CoZS}}
\newcommand{\bcozs}{\mathsf{bCoZS}}
\newcommand{\os}{\mathfrak{o}}

\newtheorem{thm}{Theorem}[section]
 \newtheorem{theorem}[thm]{Theorem}
 \newtheorem{lemma}[thm]{Lemma}
 \newtheorem{proposition}[thm]{Proposition}
 \newtheorem{corollary}[thm]{Corollary}
 \theoremstyle{definition}
 \newtheorem{definition}[thm]{Definition}
 \newtheorem{remark}[thm]{Remark}

 \numberwithin{equation}{section}
 \newtheorem{fact}[thm]{Fact}
 
\begin{document}

\title{The \v{C}ech-Stone compactification of the reals and coz-included sublocales}

\author[A. B. Avilez]{Ana Bel\'en Avilez}
\address[Avilez]{The Center for Excellence in Computation, Algebra, and Topology, Chapman University,
Orange, California 92866, USA.}
\email{avilezgarcia@chapman.edu}
\author[O. Ighedo]{Oghenetega Ighedo}
\address[Ighedo]{The Center for Excellence in Computation, Algebra, and Topology, Chapman University,
Orange, California 92866, U.S.A.}
\email{ighedo@chapman.edu}
\author[J. Walters-Wayland]{Joanne Walters-Wayland}
\address[Walters-Wayland]{The Center for Excellence in Computation, Algebra, and Topology, Chapman University,
Orange, California 92866, U.S.A.}
\email{joanne@waylands.com}
\begin{thanks}{File name:Beta and coz-inclusions}
\end{thanks}
\date{February 14, 2025}
\subjclass[2020]{18F60; 18F70; 06D22; 06E25; 54E05; 54D10}
\keywords{Topology, frame, sublocales, quotients, C-embedding, booleanization, cozero}






\begin{abstract}
		A new notion of a coz-included sublocale is introduced. It is a necessary condition  for a sublocale to be a dense $C$-embedding, and together with the property of being well-embedded, is sufficient. However, although it is clearly a strengthening of coz-ontoness (or z-embeddings), it is incomparable with $C^*$-embeddings as is demonstrated by the \v{C}ech-Stone compactification of the real numbers.  Motivated by non-trivial spatial examples, this new notion and other weaker variants, such as casi coz-included and relative zero sublocales,  are explored and characterized in the point-free setting. As a consequence, new characterizations of several classes of frames, such as perfectly normal and $Oz$ frames, among others, are provided. Besides the lattice theoretic interpretations, there are some interesting topological ramifications. 
\end{abstract}

\maketitle

\section{Introduction}
Dense $C$-embeddings have been a focus of interest in both topological spaces and their generalizations, locales (see, for example \cite{gillman_rings_1960, ball_c-_2002}. In the latter, they are closely connected to the Lindel\"of reflection (a construction that is rarely spatial): a sublocale is always densely $C$-embedded in its Lindel\"of reflection (cf. its compact coreflection in which it is densely $C^*$-embedded). In fact,     
 a sublocale  $S \subseteq L$ is a dense $C$-embedding iff $\coz S = \coz L$; this is made explicit in \cite{ball_lindelof_2017} but follows from \cite{ball_c-_2002}. 
 
 Thinking of the latter simply as a set equality, one can obviously  
split the equality into two inequalities:
$\coz L \subseteq \coz S$ and $\coz S \subseteq \coz L$.
The first obviously implies density and it captures the notion of a dense well-embedding, that is a dense, almost coz-codense frame quotient (or equivalently a dense, coz-codense quotient). See \cite{ball_well-embedding_2006}.
The second inclusion does not imply/require density and we propose to define it as $S$ being a \emph{coz-included} sublocale of $L$. By design, a dense well-embedded sublocale satisfying this condition is $C$-embedded. It has been established that $C^*$-embeddings (and even the weaker notion of coz-embeddings) together with well-embeddings are equivalent to $C$-embeddings \cite{ball_c-_2002,ball_well-embedding_2006}. This motivates understanding how ``coz-inclusions'' fit into, and help explain, the overall picture. In particular, we are able to address the errors that were made in \cite[p.~354]{ball_lindelof_2017} -- the authors state, erroneously that being a coz-embedding is equivalent to being coz-included. The subtlety of the difference between our new notion and that of coz (or $C^*$)-embeddings is highlighted by the fact that the locale of the reals is not coz-included in its \v Cech-Stone compactification (see \cref{beta R not cozincl}). One particularly pleasing consequence of this new definition is a characterization of perfect normality in terms of a condition on all sublocales (see \cref{PN iff all subs coz incl}). This motivates looking at weaker variants of perfect normality, including Oz and weak Oz frames and introducing some new ones. This allows us to correct and improve Proposition.~6.2 of \cite{ball_lindelof_2017}. \footnote{Both the errors that we correct do not effect the main results in that paper.} 

During our analysis, it became clear that the notion of being coz-embedded is useful to ensure ``good'' behavior of the cozero part of the sublocale (or quotient frame) - many of the properties of interest are only concerned with the cozeros that are inherited from the superlocale (images of cozeros of the domain of the quotient). Thus we propose to identify them as the ``bueno'' cozeros, which form a cozero base (see \cite{banaschewski_cozero_2001}) in their own right. For a sublocale $S$ of $L$, the bueno cozeros will be denoted by $\bcoz S$. This leads to a weakening of the special embeddings mentioned above. We introduce the notion of   casi coz-inclusions: $\bcoz S \subseteq \coz L$.\footnote{Weaker variants of $C^*$ and $C$-embeddings will be discussed in future work.} It also became clear that the bottom of the sublocale being cozero in the superlocale (that is, for $S$ a sublocale of $L$, $0_S \in \coz L$), was sufficient to obtain many of the results. This led to the notion of a relative zero sublocale: a sublocale whose closure is a zero sublocale (akin to relative compactness in spaces). This is also, obviously, a weakening of a sublocale being dense; the closure is the entire locale, or equivalently the bottom of the sublocale is the same as that of the superlocale (that is, $0_S= 0_L$). Noteworthy are characterizations of weaker variants of perfect normality in terms of the interaction between relative zeros and cozero elements. 

We conclude the paper with a section summarizing our results and other known results which characterize when the Booleanization  (a.k.a. skeleton) is embedded in a certain way, and then when all designated sublocales are embedded in a specified way.

\section{Preliminaries}
We assume basic knowledge of the category of frames, more specifically the full subcategory of completely regular frames, and refer the reader to, and follow the notation of, \cite{picado_frames_2012}.
Of particular importance are the following definitions/results for a frame $L$:
\begin{itemize}
\item $a \in L$ is a \emph{cozero element} if $a=\bigvee\{a_n \mid a_n \cb a\}$ (see, for example, \cite{banaschewski_real_1997}), where $\cb$ denotes the completly below relation. Let $\coz L$ denote the sub-$\sigma$-frame of cozero elements.
\item $L$  is \emph{completely regular} if $\coz L$ join generates $L$. 
    \item $L$ is \emph{normal} if whenever $ a \join b = e$, there exists $x, y$ with $x \meet y = 0$ and $x \join a = e = y \join b$. 
    \item $L$ is \emph{perfectly normal} if $L = \coz L$.
    \item $L$ is \emph{compact/Lindel\"of} if every cover of $L$ has a finite/countable subcover (and this can be extended in the obvious way to elements).
    \item The pseudocomplement of an element $a \in L$ is $a^* := \bigvee \{x \in L| a \meet x = 0\}$. 
    \item Every frame $L$ has a smallest dense sublocale, namely its Booleanization (also known as its skeleton), denoted by $L^{**}:=\{a^{**}\mid a \in L\}$.\footnote{This sublocale is often denoted by $\mathfrak B L$.} Note that, as sets, $L^{**}= L^*$.
\end{itemize}
\begin{remark} The regular $\sigma$-frame $\coz L$ is always normal (since regularity implies normality for $\sigma$-frames), and thus perfect normality implies normality for completely regular frames.
    
\end{remark}

It is often useful and enlightening to use the alternate descriptions of sublocales, namely as quotients, nuclei, or congruences. And we freely use the different descriptions to prove results from various perspectives -- sometimes because these are often easier or more illustrative.
We include the details below to establish notation and remind the reader of some of the essential relationships:
\begin{itemize}
    \item Given a frame map $h\colon L\to M$ we denote its right adjoint by the localic map $h_*\colon M\to L$, and we write $h\dashv h_*$.  
    \item If $S$ is a sublocale of $L$ with the localic embedding $j\colon S \hookrightarrow L$, then $\map{q}{L}{S}$ denotes the associated frame quotient (meaning $j=q_*$ and $q\dashv j$), $jq$ is the corresponding nucleus, and $$E_S=E_{q} = \{(x,y):q(x)=q(y)\}$$ is a congruence. In particular, $x \in S$ iff $x = \bigvee E_S[x]$.
\item If $\map{h}{L}{M}$ is a frame quotient, then $h_*[M]$ is a sublocale of $L$, $h_*h$ is a nucleus, and $E_h$ is a congruence. Note that $h_*[M]\cong M$ and the isomorphism is given by $x\mapsto h(x)$ for every $x\in h_*[M]$. Consequently, $\coz M\cong \coz h_*[M] $.

\end{itemize}

\begin{lemma} For any quotient $\map{h}{L}{M}$ and associated sublocale $S$ of $L$ 
\label{cozsh}
    \[
        \coz S =\coz h_*[M]= h_*[\coz M].
   \]
        
\end{lemma}

\begin{proof}
    Since $S= h_*h[L]$ and $h[L]=M$, the first equality is clear. If $x \in \coz S$, then $x \in S$ so $x = h_*h(x)$. But $h(x) \in \coz M$, since $S\cong M$ and the isomorphism is given by $h$. Thus $x=h_*h(x) \in h_*[\coz M]$. For the other inclusion, if $x \in h_*[\coz M]$, then $x=h_*(y)$ for some $y \in \coz M$. The map  $h_*:M \longrightarrow h_*[M] =S$ is a localic isomorphism because $h$ is a frame quotient. Since adjoints are unique,  $h_*$ is also a frame isomorphism and thus preserves cozeros. Hence, $h_* (y) \in \coz h_*[M]=\coz S$, so $x \in \coz S$. 
\end{proof}
\begin{lemma}\label{coz q*} For any frame map $\map{h}{L}{M}$, if $\map{\coz h}{\coz L}{\coz M}$ has a right adjoint, then it is $h_*$ restricted to $\coz M$.
    
\end{lemma}
\begin{proof} Suppose $\map{f}{\coz M}{\coz L}$ is the right adjoint to $\coz h$. We claim that $f(b) = h_* (b)$ for all $b \in \coz M$:
for any $a\in \coz L$, 
  \[a\leq f(b)\iff \coz h(a)\leq b \iff h(a)\leq b \iff a\leq h_*(b).\]
  So by complete regularity of $L$, we have $f(b)=h_*(b)$. Thus, $f= h_*|_{\coz M} $.
    
\end{proof}
We caution the reader that this adjoint does not necessarily exist: the restriction of $h_*$ to $\coz M$ always exists but need not map into $\coz L$. See \cref{casi Coz incl charact.} below. When it does exist, we will denote the right adjoint by $\coz h_*$.

\begin{lemma}\label{coz nucleus} For any nucleus $\map{\upsilon}{L}{L}$,\\
$$\upsilon \text{ drops to a nucleus on } \coz L \iff  \upsilon (a) \in \coz L \text{ for each } a \in \coz L.$$
\end{lemma}
\begin{proof} The restriction of $\upsilon$ to $\coz L$, $\map{\upsilon |_{\coz L}}{\coz L}{L}$ always satisfies the conditions of being a nucleus except that the codomain need not be $\coz L$.
\end{proof}
When a nucleus $\upsilon$ drops as above, we will denote it by $\coz \upsilon$.\\

We remind the reader of some basic notation for sublocales \cite[III.6]{picado_frames_2012}:
for any $a \in L$, $\close_L(a)$ and $\open_L(a)$ denote the closed and open sublocales of $L$.\footnote{We shall denote simply by $\close(a)$ and $\open(a)$ when there is no danger of confusion.} If $a \in \coz L$, these are referred to as zero or cozero sublocales. The collections of closed, open, zero, and cozero sublocales of $L$ are denoted by 
$$\close(L),\quad \open(L),\quad \zs(L),\quad \text{and}\quad\cozs(L)\quad \text{respectively.}$$
So, in particular, $$\close_L(a)\in \zs (L) \iff \open_L(a)\in \cozs(L)\iff a\in \coz L.$$
The closure $\overline S$ of a sublocale $S$ is the smallest closed sublocale containing $S$, given by the formula $\overline S = {\uparrow} (\bigwedge S)$.
The interior $S^\circ$ is the largest open sublocale contained in $S$.  In particular, $\overline{\open(a)}=\close(a^*)$ and, consequently, $\close(a)^\circ=\open(a^*)$.

We also need the following definitions and/or characterizations of sublocales: a sublocale $S$ of $L$ with associated quotient map $\map{q}{L}{S}$ is a 
\begin{itemize}
    \item coz-embedding if $q$ is coz-onto (see \cite[p.~6]{ball_c-_2002}); \footnote{These are also called ``z-embeddings" (see, for e.g., \cite[Sec.~6]{gutierrez_garcia_notes_2014}).}
    \item $C^*$-embedding if every bounded continuous real function to $S$ extends to $L$ (see \cite[Defn.~7.0.11]{ball_c-_2002}); 
    \item $C$-embedding if every continuous real function to $S$ extends to $L$ (see \cite[Defn.~7.0.11]{ball_c-_2002}). 
    
\end{itemize}
This is analogous to the situation in spaces.
\begin{fact}\label{char. of normality}
    \begin{align*}
 \mbox{A frame is normal}\; &\iff \text{all closed sublocales are coz-embedded \cite[4.11]{dube_coz-onto_2006}};\\
  &\iff \text{all closed sublocales are C*-embedded \cite[8.3.3]{ball_c-_2002}};\\
    &\iff \text{all closed sublocales are C-embedded \cite[8.3.3]{ball_c-_2002}.}\\
\end{align*}
\end{fact}
Compact regular frames form a coreflective subcategory of all completely regular frames (see \cite[VII.4]{picado_frames_2012}). The coreflection $\beta L$ (a.k.a. the \v{C}ech-Stone compactification) may be described as the frame of regular ideals: an ideal $J$ of $L$ is \emph{regular} if for each $x\in J$, there exists $y \in J$ with $x\cb y$.\footnote{These ideals are often called completely regular or round ideals.} The coreflection map is then given by join, and we denote its right adjoint by $\DDarrow$ where ${\DDarrow} a := \{x\in L\mid x \cb a\}$. Thus, $L$ is a quotient of $\beta L$ and so may be identified as a sublocale of $\beta L$. In particular, it is a dense $C^*$-embedded (and hence coz-embedded) sublocale (see, for example, \cite[Cor. 8.2.7]{ball_c-_2002}). Moreover, they have the same Booleanization. Cozeros of $\beta L$ are precisely the countably generated regular ideals (see \cite{banaschewski_stonecech_1989}). 
In $\beta L$, $I \prec J$ iff $I \subseteq {\DDarrow} u$ for some $u \in J$ iff $I \subseteq {\DDarrow} u \subseteq J$ for some $u \in L$. Compact regular frames have many desirable properties; in particular they are normal.

\section{The locale of the reals and its \v{C}ech-Stone compactification}

 The frame of the reals, denoted $\mathfrak L \Rn$, is presented in detail in \cite[XIV]{picado_frames_2012} and \cite{banaschewski_real_1997}. It is worth noting that $\mathfrak L \Rn$ is isomorphic to the usual topology of the reals.  Of particular importance in this paper are the following properties of $\mathfrak L \Rn$:
\begin{itemize}
    \item $\mathfrak L \Rn$ is completely regular;
    \item $\mathfrak L \Rn$ is Lindel\"of;
    \item $\mathfrak L \Rn$ is perfectly normal;
    \item $\mathfrak L \Rn \cong \mathfrak{O} \Rn$;
    \item the Booleanization, $ (\mathfrak L \Rn)^{**}$, consists of all regular open sets.
\end{itemize}

We can identify $\mathfrak L \Rn$ as a dense $C^*$-embedded sublocale of its \v{C}ech-Stone compactification, $\beta \mathfrak L \Rn$,\footnote{$\beta \mathfrak L \Rn \cong \mathcal{O}\beta \Rn \cong \beta \mathcal O \Rn $
(with appropriate choice principles).} and their shared Booleanization may be identified  as $\{{\DDarrow} a\mid a \in (\mathfrak L \Rn)^{**}\}$.

\begin{proposition} $\beta \mathfrak L \Rn$ is not perfectly normal.\label{beta R not PN}
    
\end{proposition}
\begin{proof}

We will demonstrate this by showing the existence of a regular ideal that is not cozero.
    Consider the open set $v= \bigvee_{m \in \N} (2m,2m+1)$.
We claim that the regular ideal ${\DDarrow} v$ is not cozero (that is, not countably generated): 

Let $I_n$ be any countable sequence of regular ideals completely below ${\DDarrow} v$. W.l.o.g. we can arrange things so that
\begin{enumerate}
    \item each $I_n = \,{\DDarrow} u_n$ for some $u_n$, 
    \item $u_n \prec u_{n+1}$ (if $x, y \prec a$ then $x\vee y \prec a$).
\end{enumerate}
We use the following facts:
\begin{itemize}
    \item in $\mathfrak{L} \Rn$ , $(p,q) \prec (r,s)$ iff $r<p<q<s$.
    \item  $v$ is a union of properly disjoint intervals, and thus so is each $u_n$.
\end{itemize}
 Since $u_n \prec v$, $u_n \meet (2m,2m+1) \prec (2m,2m+1)$. So
for each $m$, let \[w_{m,n} = (2m,2m+1) \meet u_n\] 

Then it follows that $w_{m,n} \prec w_{m,n+1} \prec (2m,2m+1)$.

Consider the element $w=\bigvee_k w_{k,k} $.
\begin{itemize}
    \item  $w \prec v$ : $ v^*  \vee \bigvee ((w_{k,k})^* \meet (2k-\frac{1}{4}, 2k+\frac{5}{4}))$
 is a separating element.
 \item $w$ is not below any $u_n$:
for any $n$, if $k >n$ then $w_{k,k} > (2k,2k+1) \meet u_n$
\end{itemize}
Thus $w \in {\DDarrow} v$, but $w \notin \bigvee I_n$, proving that ${\DDarrow}v$ is not cozero.
\end{proof}
\begin{remark}\label{beta R not cozincl} Considering $\mathfrak L \Rn$ as a sublocale of $ \beta \mathfrak L \Rn$ (that is, $\mathfrak L \Rn\subseteq \beta \mathfrak L \Rn$), the proof above shows that $\mathfrak L \Rn$ is not a subset of $\coz \beta \mathfrak L \Rn$. Moreover, since $\mathfrak L \Rn$ is perfectly normal, $\mathfrak L \Rn = \coz(\mathfrak L \Rn)$, so $\coz(\mathfrak L \Rn) \not\subseteq\coz \beta \mathfrak L \Rn $.
    
\end{remark}
    
\begin{definition}\label{defn: coole}
    A frame is \emph{coole} if its Booleanization is coz-embedded; equivalently, the booleanization map ${L} \rightarrow{L^{**}}$ is coz-onto.
\end{definition}

\begin{corollary} \label{betaR is coole} $\beta \mathfrak L \Rn$ is a coole frame.
    
\end{corollary}

\begin{proof} Recall that $(\mathfrak L \Rn)^{**}=(\beta \mathfrak L \Rn)^{**}$, so the result follows since both $\beta \mathfrak L \Rn \to \mathfrak L \Rn $ and $\mathfrak L \Rn \to (\mathfrak L \Rn)^{**}$ are coz-onto. The former is always true for any \v{C}ech-Stone compactification, the latter since $\mathfrak L \Rn $ is perfectly normal.
\end{proof}
\begin{remark}\label{ex. coole not Oz}The element ${\DDarrow}v$
from the proof of \cref{beta R not PN} is a regular element but not a cozero of $\beta \mathfrak L \Rn $; that is, $(\mathfrak L \Rn)^{**} \not\subseteq\coz \beta \mathfrak L \Rn $. 
    
\end{remark}

\section{Coz-inclusions}\label{sec coz inclusions}

As mentioned in the introduction, an obvious weakening of being a dense $C$-embedded sublocale is the following:

\begin{definition} A sublocale $S \subseteq L$ is \emph{coz-included} if, as sets, 
$\coz S \subseteq \coz L$.  
\end{definition}
We abuse terminology and say a frame quotient $\map{h}{L}{M}$ is coz-included if the associated sublocale $h_*[M]$ is coz-included, that is, ${\coz h_*[M]\subseteq \coz L}$.
 
\begin{proposition}\label{cozincl implies cozemb} Coz-inclusions are coz-embeddings but not conversely.
\end{proposition}
\begin{proof} Suppose $\coz S \subseteq \coz L$ and take $x \in \coz S$. Then $x \in \coz L$, by assumption. Moreover, since $x \in S$, $x=q(x)$ and thus $q$ is coz-onto. For the negation of the converse, consider $\mathfrak{L} \Rn$ as a sublocale of $\beta(\mathfrak{L} \Rn)$: this is coz-embedded but not coz-included; see \cref{beta R not cozincl}.
\end{proof}
\begin{remark}\
 We point out that in \cite[Pg.~354]{ball_lindelof_2017}, the authors erroneously claimed that the two conditions were equivalent. This led to an incorrect conclusion in \cite[Prop.~6.2]{ball_lindelof_2017} which will be remedied below (see \cref{char of Oz}).

\end{remark}

Decoding the definition above gives the following characterization.
\begin{lemma} \label{coz incl iff adj cozeros}A quotient $\map{h}{L}{M}$ is coz-included iff $h_*$ preserves cozeros.
    
\end{lemma}
\begin{proof} By \cref{cozsh}, $ h_*[\coz M]= \coz h_*[M]$. Since coz-included means $\coz h_*[M]\subseteq \coz L$, and $h_*$ preserves cozeros means $h_*[\coz M]\subseteq \coz L$, the two notions are clearly equivalent.
\end{proof}
We now give various equivalent conditions, including what coz-included means for the associated nucleus and congruence.
\begin{theorem}\label{Coz incl charact.} For any sublocale $S$ of $L$, the following are equivalent:
 \begin{enumerate}[label=\normalfont(\arabic*)]
     \item \label{CI:part1} $S$ is coz-included.
     \item \label{CI:part2} For every zero sublocale $Z$ of $S$, the closure in $L$ of $Z$ is a zero sublocale of $L$; that is, $Z\in \zs(S)$ implies $\overline{Z}^L \in \zs(L)$.
     \item \label{CI:part3} The quotient $\map{q}{L}{S}$ is coz-included.
     \item\label{CI:part4} The right adjoint, $q_*$, of the quotient map preserves cozeros.
      \item\label{CI:part5} The map $\map{\coz q}{\coz L}{\coz S}$ has a right adjoint.
      \item \label{CI:part6} $S$ is coz-embedded and the associated nucleus $\upsilon_S$ drops to a nucleus on $\coz L$ (namely $\coz \upsilon_S$). 
      \item\label{CI:part7} $c \in  \coz S$ implies $\bigvee E_S[c] \in \coz L$.
 \end{enumerate}    
 \end{theorem}
 \begin{proof}

\ref{CI:part1} $\implies$ \ref{CI:part2}: Let $Z$ be a zero sublocale of $S$; that is, $Z=\close_S(a)$ for some $a\in \coz S$. Then $\overline{Z}^L=\overline{\close_S(a)}^L=\close_L(a)$.
    By assumption $a\in \coz L$, so $\close_L(a)=\overline{Z^L}$ is in $\zs(L)$.\\
    \ref{CI:part2}$ \implies $ \ref{CI:part1}: Let $a\in \coz S$, then $\close_S(a)$ is a zero sublocale of $S$. Hence, by assumption, $\overline{\close_S(a)}^L$ is a zero sublocale of $L$. But $\overline{\close_S(a)}^L=\close_L(a)$, yielding $a \in \coz L$.\\
 \ref{CI:part1}$\iff$ \ref{CI:part3}: By definition.\\
 \ref{CI:part3}$\iff$\ref{CI:part4}: This is precisely  \cref{coz incl iff adj cozeros}.\\
 \ref{CI:part4} $\implies$\ref{CI:part5}: If $q_*$ preserves cozero elements, then the restriction $q_*|_{\coz S}$ maps to $\coz L$. Furthermore, since $q\dashv q_*$, for every $a\in \coz L$ and $b\in \coz S$ we have that 
    \[ \coz q (a)= q(a)\leq b\iff a\leq q_*(b)=q_*|_{\coz S} (b).\]
Thus $\coz q$ has a right adjoint namely, $\coz q_*$ (see \cref{coz q*}).\\
 \ref{CI:part5}$\implies$ \ref{CI:part6}: 
 By \cref{coz q*}, $\map{\coz q_*}{\coz S}{\coz L}$ is the right adjoint to $\coz q$. Take $x \in \coz S$, then $\coz q_*(x) \in \coz L$, by assumption. But $q(\coz q_* (x))= qq_*(x) = x$, showing that $S$ is coz-embedded. By \cref{coz nucleus}, $\upsilon_S$ drops to a nucleus on $\coz L$ iff for each $y\in \coz L$, $\upsilon_S(y)\in \coz L$. Suppose $y \in \coz L$, then $\upsilon_S(y)=q_*(q(y))=\coz q_*(q(y))$ (the last equality holds since $q(y)\in \coz S$). But, by assumption $\map{\coz q_*}{\coz S}{\coz L}$ and thus $\upsilon_S(y) \in \coz L$. 
     
 \ref{CI:part6}$\implies$\ref{CI:part4}: Take $x \in \coz S$, then $x = q(y)$ for some $y \in \coz L$ since $S$ is coz-embedded in $L$. Now $q_*(x) = q_*(q(y))= \upsilon_S(y) = \coz \upsilon_S(y)$ since $y \in \coz L$, and thus by the assumption, this element is in $\coz L$, as desired.

  \ref{CI:part7}$\iff$\ref{CI:part1}: Note that for every $c\in \coz S$, since $c\in S$  we have that $c=\tbigvee E_S[c]$. Thus, $\tbigvee E_S[c]\in \coz L$ for every $c\in \coz S$ if and only if $\coz S\subseteq \coz L$. 
 \end{proof}
 
  \begin{remark}\
  \begin{enumerate}[label=\normalfont(\arabic*), ref = \theproposition(\arabic*)]
  \item   Recall that $S$ being coz-embedded means every zero sublocale $Z$ of $S$ is the intersection of $S$ with a zero sublocale of $L$, whereas coz-included is that for every zero sublocale $Z$ of $S$ there is a smallest zero sublocale $Z'$ of $L$ such that $Z'\cap S=Z$. There is an analogous comparison for cozero elements: $S$ is coz-embedded if every cozero of $S$ is the image of a cozero of $L$, whereas coz-included means that the largest element mapped to a given cozero of $S$ is a cozero of $L$.
  \item For a topological interpretation, a $z$-embedded subspace $S$ of a topological space $X$ is a subset where every zero set is an intersection of $S$ with a zero set of $X$. This is the desirable, but not always common  behavior of zero sets of a subspace. If for every zero set of $S$ we can find the smallest zero set of $X$ (namely, its closure in the in $X$) that satisfies this property, we say the subspace is coz-included.
      \item  If $S$ is coz-included, then, since it is obviously a zero sublocale of itself, $\overline S$ is a zero sublocale of $L$. Obviously this condition is trivially true for any dense sublocale, but this is not the case for every coz/$C^*$-embedded sublocale. In other words, the closure of a coz/$C^*$-embedded sublocale may not be a zero sublocale. For instance consider a normal, not perfectly normal frame (for example, $\beta \mathfrak{L} \Rn$), and take a proper closed sublocale that is not a zero sublocale.
       \end{enumerate}
  \end{remark}
  \begin{proposition}\label{when zeros are ci} For a completelely regular frame $L$,
\begin{enumerate}[label=\normalfont(\arabic*), ref = \theproposition(\arabic*)]
    \item \cite[3.2.11]{ball_c-_2002} every cozero sublocale is coz-embedded;
    \item a zero sublocale is coz-embedded iff it coz-included.\label[proposition]{when zeros are ci 2}
\end{enumerate}
    
\end{proposition}
\begin{proof}
    (2) Coz-included always implies coz-embedded (\cref{cozincl implies cozemb}).
    Conversely, let $Z$ be a zero sublocale of $L$, and assume it is coz-embedded. To show it is coz-included, let $F$ be a zero sublocale of $Z$. Since $Z$ is coz-embedded in $L$, $F=Z'\cap Z$ for some zero sublocale $Z'$ of $L$. Thus $F$ is a zero sublocale of $L$ and $F=\overline{F}$, and thus \ref{CI:part2} of \cref{Coz incl charact.} holds.
    \end{proof}

Recall that a frame is 
\emph{weakly perfectly normal} if all sublocales are coz-embedded \cite[P.74]{avilez_point-free_2022}.
   The latter is the obvious generalization of the spatial notion introduced in \cite{blair_extensions_1974}, and it is indeed a weaker notion than that of perfect normality \cite[Prop.~4.3.6]{avilez_point-free_2022}.
We obtain the following rather pleasing characterization of perfect normality which helps explain this difference.
\begin{theorem}\label{PN iff all subs coz incl} A completely regular frame is perfectly normal iff every sublocale is coz-included iff every closed sublocale is coz-included.
    
\end{theorem}
\begin{proof} $L$ is perfectly normal means $L = \coz L$: any sublocale $S$ is clearly a subset of $L$, and thus so is the cozero part yielding that each sublocale is coz-included. Now take any $a\in L$, and consider $\close_L(a)$. Clearly $a \in \coz \close_L(a)$ (since it is the bottom of that sublocale), and hence $a \in \coz L$ since $\close_L(a)$ is coz -included in $L$, by assumption.
\end{proof}
\begin{remark}\label{skel of reals not coz-incl}\
 The Booleanization of any perfectly normal frame is coz-included, specifically so is that of the reals. However, as was shown \cref{ex. coole not Oz}, the Booleanization of the reals is not coz-included in $\beta \mathfrak L \Rn$.
\end{remark}
The characterization in \cref{PN iff all subs coz incl}
 also emphasizes the difference with normality, see \cref{char. of normality}
and leads us to consider how coz-inclusions relate to $C$ and $C^*$-embeddings. 
     We remind the reader of the following, where $CL$ denotes the (Boolean) lattice of complemented elements of $L$:
     a frame $L$ is \emph{extremally disconnected}  if $CL = L^{**}$. Moreover, a frame is extremally disconnected iff all open sublocales are $C^*$-embedded iff all dense sublocales are $C^*$-embedded \cite[8.4.1]{ball_c-_2002}.   In particular, a frame is extremally disconnected iff its Booleanization is $C^*$-embedded.   

\begin{proposition} \label{coz incl vs Cstar}
Coz-inclusions are not necessarily $C^*$-embeddings, nor conversely.

\end{proposition}
\begin{proof}

$\mathfrak{L}\Rn \subset \beta(\mathfrak{L} \Rn)$ is $C^*$-embedded but not coz-included by \cref{beta R not cozincl}.
         $\mathfrak{L} \Rn$ is perfectly normal but not extremally disconnected, so there exists an open (in fact, any proper open interval) that is not $C^*$-embedded but obviously coz-included by \cref{PN iff all subs coz incl}.
 \end{proof}
\begin{proposition}\label{dense C implies CI}
    Dense $C$-embeddings are coz-inclusions, but C-embeddings are not necessarily so.
\end{proposition}
\begin{proof}
    If a dense $S$ is C-embedded in $L$, then $\coz S = \coz L$ by \cite[Thm. 5.5]{ball_lindelof_2017}, and thus, in particular $S$ is coz-included in $L$. A normal frame that is not perfectly normal (for example, $\beta \mathfrak L \Rn$) will have a closed sublocale that is C-embedded but not coz-included by \cref{PN iff all subs coz incl}. 
\end{proof}

\cref{fig:coz incl vs emb vs C*} summarizes the above implications. Note that these are all strict with differentiating examples given below. Moreover, density of the sublocale is required to ensure the lower left implication. Indeed coz-inclusions and C-embeddings are not comparable without density (see \cref{dense C implies CI}).

\begin{figure}[h]
  \[\begin{tikzcd}
	& \text{coz-embeddings} \\
    \text{coz-inclusions} && \text{$C^*$-embeddings} \\
    & \text{dense $C$-embeddings}
	\arrow["\diamondsuit",Rightarrow, from=2-1, to=1-2]
	\arrow["\heartsuit"',Rightarrow, from=2-3, to=1-2]
	\arrow["\clubsuit", Rightarrow, from=3-2, to=2-1]
	\arrow["\spadesuit"',Rightarrow, from=3-2, to=2-3]
\end{tikzcd}\]
\caption{}
    \label{fig:coz incl vs emb vs C*}
\end{figure}

\begin{enumerate}
    \item[] $\spadesuit$: $\mathfrak L \Rn \subseteq \beta \mathfrak L \Rn$ is $C^*$-embedded, but not C-embedded.
    \item[] $\diamondsuit$: $\mathfrak L\Rn \subseteq \beta \mathfrak L\Rn$ is coz-embedded, but not coz-included.
    \item[] $\heartsuit$: every proper open sublocale of $\mathfrak L \Rn$ is coz-included, but not $C^*$-embedded. 
    \item[] $\clubsuit$: the Booleanization of $\mathfrak L \Rn$ is coz-included, but not C-embedded (since $\mathfrak L \Rn$ is perfectly normal but not extremally disconnected). 
\end{enumerate}

\section{The locale of the one-point compactification of an uncountable set}
 \label{1-pt comp}
 Consider an uncountable discrete set $E$ and its one-point compactification $X = E \cup \{\infty\}$.  $\O X$, the frame of open subsets of $X$, is the family of all subsets $C \subseteq X$ such that if $\infty \in C$ then $C$ must be cofinite. This situation is considered in \cite[Example 3.1]{blair_spaces_1976}).
 
 The associated frame quotient  $h:\O X \rightarrow \O E$  is defined by
$h(A) = A \cap E$,
and hence 
\[h_*(B) =\begin{cases}
    B &\mbox{  if } B \ \mbox {is not cofinite}\\
    B\cup \{\infty\}& \mbox{otherwise.}
\end{cases} \] 
 Notice that the frame $\O E$ is not a sublocale of $\O X$, but it is isomorphic to $\{B\subset E: B \text{ not cofinite}\}\cup\{\{B\cup\{\infty\}\}: B \text{ cofinite} \}$ which is a sublocale of $\O X$. We identify $\O E$ with this sublocale.\\
     
\noindent{ \bf The cozero sets of $\O E$ and $\O X$}: $E$ is discrete so $\coz \O E = \O E$. Since $\O X$ is compact, the cozeros are precisely the Lindel\"of elements \footnote{In a regular Lindel\"of frame, the cozero elements are precisely the Lindel\"of ones.}. Thus the cozeros are exactly the cofinite sets containing $\infty$ together with the countable sets which omit $\infty$. 

 In \cite{blair_spaces_1976}, they show the embedding is not a coz-embedding,\footnote{We will obtain a ``frame" version of this in \cref{coz incl iff casi and Z}.} and thus by \cref{cozincl implies cozemb} it is not coz-included. For the sake of completeness, and as it is an essential concept in the development and justification of the next section, we give the following direct proof.
 \begin{proposition} \label{1pt comp not cozincl}
For $E$ and $X$ as described above, $\O E$ is not coz-included in $\O X$.
\end{proposition}
    
     \begin{proof} By \cref{coz incl iff adj cozeros} it suffices to show that $h_*$ does not preserve cozeros. Let $A\subseteq E$ be a cocountable, not cofinite subset of $E$, then $h_*(A) = A$. Clearly $A \in \coz \O E$ but, since $A$ is not countable (and $\infty \notin A$), $A$ is not a cozero of $\O X$. Thus,  $\coz \O E \nsubseteq \coz \O X$.
     \end{proof}

\noindent{\bf The Booleanization of $\O E$ and $\O X$}: Since $E$ is discrete, $\O E = (\O E)^{**}$. Moreover, $(\O X)^{**} = (\O E)^{**}$ since $(\O E)$ is dense in $(\O X)$. Hence $(\O X)^{**} = \O E$.
In fact,  
for $A \in \O X$, 
\[A^{**} = \begin{cases}
    A &\mbox{if } A \mbox{ is not cofinite,}\\
    A\cup\{\infty\} &\mbox{if } A \ \mbox{is cofinite}.
\end{cases}
\]

\begin{proposition}\label{prop: cozincl vs casi}\
  For $X$ as described above,  $(\coz \O X)^{**} \subseteq \coz \O X$. 
\end{proposition}
 \begin{proof} Since the cozeros are precisely the cofinite sets containing $\infty$ together with the countable sets which omit $\infty$. $A^{**}=A$ for cozero $A$, and hence the result.
 \end{proof}
 
      \begin{remark}
          Note that $h$ and the Booleanization map $(A\mapsto A^{**})$ determine the same sublocale of $\mathcal{O}X$ which is isomorphic to $\mathcal{O}E$. The booleanization map is a frame quotient so it preserves cozero elements. However, not every cozero of the Booleanization (recall that $\coz (\mathcal{O}X)^{**}=(\mathcal{O}X)^{**}$ because the Booleanization is a boolean frame) is an image of a cozero. For instance, take the same $A\subseteq E$ as in the proof of \cref{1pt comp not cozincl}, that is a cocountable but not cofinite subset of $E$. It is a (cozero) element of the Booleanization ($A=A^{**}$) but one can readily see that it cannot be the image of a cozero element of $\mathcal{O}X$ because these are the countable sets that omit $\infty$ or the cofinite sets toegether with $\infty$. This means the Booleanization map is not surjective on its cozero part; equivalently $(\mathcal{O}X)^{**}$ is not coz-embedded in $\mathcal{O}X$.\color{black}
          \end{remark}
\section{Casi coz-inclusions}
The above example makes it clear that there should be a way to think about ``coz-inclusions" that are not coz-embeddings; that is, to find the essence of being a coz-inclusion. To that end, we need the notion of a ``good" cozero of a sublocale, that is, one that is realized by a cozero from the locale (and hence, every cozero is ``good" iff the sublocale is coz-embedded):

 \begin{definition}\
 \begin{enumerate}
     \item For a frame map $\map{h}{L}{M}$,
    an element $x$ of $M$ is an \emph{h-bueno cozero} if {$x~ \in~h[\coz L]$}.
    \item For a sublocale $S \subseteq L$, $x \in q[\coz L]$ is simply called a \emph{bueno cozero} of $S$, where $q$ is the associated quotient to the localic embedding of $S$ in $L$. Let $\bcoz S$ denote the collection of all bueno cozeros of $S$. 
 \end{enumerate} 
\end{definition}

\begin{remark}\label{rem:casi CI}\
\begin{enumerate}[label=\upshape(\arabic*), ref = \theremark(\arabic*)]
    \item Clearly $\bcoz S=q[\coz L]$ (or $\bcoz h_*[M]\cong h[\coz L] $) and $\bcoz S \subseteq \coz S$ (or $\bcoz h_*[M] \subseteq \coz h_*[M]$).
 Moreover, $\bcoz S$ is always a regular sub-$\sigma$-frame of $\coz S$ that join generates $S$, that is, a cozero base of $S$ (see \cite[Defn.~1]{banaschewski_cozero_2001}). It turns out that many notions dealing with cozeros of the sublocale may be weakened by considering the notions relative to this, possibly smaller, cozero base.
\item Since every Lindel\"of sublocale is coz-embedded \cite[Prop.~3.2]{dube_coz-onto_2006}, all cozeros of a Lindel\"of sublocale are bueno cozeros; that is, $\bcoz S = \coz S $.
\end{enumerate}
\end{remark}
There is an obvious way to extend the idea of a bueno element to sublocales:
\begin{definition}For a sublocale $S$ of $L$,
    
\label[remark]{rem:casi CI 2}
\begin{enumerate}
    \item A zero sublocale $Z$ of $S$ is \emph{bueno} if $Z=Z'\cap S$ for some zero sublocale $Z'$ of $L$. Let $\bzs (S)$ denote the collection of bueno zero sublocales of $S$. 
\item  A cozero sublocale $C$ of $S$ is \emph{bueno} if $C=C'\cap S$ for some cozero sublocale $C'$ of $L$. Let $\bcozs (S)$ denote the collection of bueno cozero sublocales of $S$.
 \end{enumerate}
 \end{definition}
 \begin{remark}\
     \begin{enumerate}
         \item $\bzs (S) = \left\{Z'\cap S\mid Z'\in\zs(L)\right\}$ is a subset of $\zs(S)$.
         \item $\bcozs (S)= \left\{C'\cap S\mid C'\in\bcozs(L)\right\}$ is a subset of $\cozs(S)$.
         \item For $a\in S$, the sublocale $\close_S(a)$ is a bueno zero sublocale of $S$ iff $\open_S(a)$ is a bueno cozero sublocale of $S$ iff $a\in \bcoz S$.
         That is, $$\close_S(a)\in \bzs (S) \iff \open_S(a)\in \bcozs(S)\iff a\in \bcoz S.$$
     \end{enumerate}
 \end{remark}

  We are now able to define a condition on a sublocale that is weaker than being coz-included, and as we will see it is precisely this condition  that is required to obtain the equivalence of coz-embeddings and coz-inclusions (see \cref{fig:coz incl vs emb vs C*}), which justifies using the spanish adverb ``casi'' \footnote{Our intent for future work is to always use the word ``casi'' (which means almost in Spanish) to weaken a notion so that together with coz-embeddings is equivalent to the notion.}. 
\begin{definition} A sublocale $S \subseteq L$ is \emph{casi coz-included} if, as sets, 
$\bcoz S \subseteq \coz L$.  
\end{definition}

  \begin{proposition}\label{coz incl iff casi and Z} $S$ is coz-included iff it is casi coz-included and coz-embedded.
      
  \end{proposition}
  \begin{proof} Obviously, coz-included implies both  casi coz-included and coz-embedded (for the latter see \cref{cozincl implies cozemb}). Suppose $S$ is casi coz-included, that is, $\bcoz S \subseteq \coz L$, and $S$ is coz-embedded, that is, $q[\coz L] = \coz S$. The result follows since $\bcoz S = q[\coz L]$.
     \end{proof}
\begin{remark}\label{example of casi not coz included}\
    \begin{enumerate}[label=\upshape(\arabic*), ref = \theremark(\arabic*)]
        \item The one point compactification of an uncountable discrete space discussed in \cref{1-pt comp} reveals the subtle difference between them: by \cref{1pt comp not cozincl},  $\O E$ is not coz-included in $\O X$, but by \cref{prop: cozincl vs casi},  $\O E$ is casi coz-included in $\O X$. (And hence, by \cref{coz incl iff casi and Z}, $\O E$ is not coz-embedded).
   \item The Booleanization of $\beta \mathfrak L\Rn $, namely $(\mathfrak L\Rn)^{**}$, is coz-embedded, but not coz-included (see \cref{ex. coole not Oz}) and hence it is NOT casi coz-included by \cref{coz incl iff casi and Z}.\label[remark]{example of not casi but coz-embedded}
 \item For Lindel\"of sublocales,  since all cozeros are bueno cozeros, casi coz-included is equivalent to coz-included.
  \end{enumerate}
\end{remark}
Again, we abuse terminology and say a quotient $\map{h}{L}{M}$ is casi coz-included if ${h_*[h[\coz L]]\subseteq \coz L}$.
\begin{lemma}\label{casi coz incl iff adj cozeros}A quotient $\map{h}{L}{M}$ is casi coz-included iff $h_*$ preserves h-bueno cozeros; that is, $h_*(a)\in \coz L$ for every $a\in h[\coz L]$.  
\end{lemma}
\begin{proof}
   Suppose $h\colon L\to M$ is casi coz-included, that is, $h_*[h[\coz L]]\subseteq \coz L$. If $a\in h[\coz L]$, then $h_*(a)\in h_*[h[\coz L]]\subseteq \coz L$, as required.\\
    Now, assume $h_*$ preserves h-bueno cozeros, and let $a\in h_*[h[\coz L]]$. Then $a=h_*h(x)$ for some $x\in \coz L$. But $h(x) \in h[\coz L]$, and thus, by assumption $h_*(h(x)) \in \coz L$ and thus so is $a$.
\end{proof}

\begin{lemma}\label{qcoz q*} For any frame map $\map{h}{L}{M}$, if $\map{\coz h}{\coz L}{h[\coz L]}$ has a right adjoint, then it is $h_*$ restricted to $h[\coz L]$.
    
\end{lemma}
\begin{proof} The same proof as in \cref{coz q*} works when restricting to bueno cozeros.
\end{proof}
As we will see below, with some subtle adjustments, the proofs in \cref{Coz incl charact.} all work with the restriction to the bueno cozeros, and thus we obtain the following characterizations of casi coz-inclusions. 
\begin{theorem}\label{casi Coz incl charact.} For any sublocale $S$ of $L$, the following are equivalent:
 \begin{enumerate}[label=\upshape(\arabic*)]
     \item \label{casi CI1}$S$ is casi coz-included.
     \item \label{casi CI2} For every bueno zero sublocale $Z$ of $S$, the closure in $L$ of $Z$ is a zero sublocale of $L$; that is, $Z \in \bzs (S)$  implies $\overline{Z}^L \in \zs(L)$.
     \label[proposition]{casi Coz incl charact.2}
     \item \label{casi CI3}The quotient $\map{q}{L}{S}$ is casi coz-included.
     \item \label{casi CI4}The right adjoint, $q_*$, of the quotient map preserves bueno cozeros.
      \item\label{casi CI5} The map $\map{\coz q}{\coz L}{\bcoz S}$ has a right adjoint.
      \item\label{casi CI6} The associated nucleus $\upsilon_S$ drops to a nucleus on $\coz L$.
      \item\label{casi CI7} $c \in  \bcoz S$ implies $\bigvee E_S[c] \in \coz L$.
       \item\label{casi CI8} $c \in  \coz L$ implies $\bigvee E_S[c] \in \coz L$.
 \end{enumerate}
     
 \end{theorem}
  \begin{proof}
\ref{casi CI1}$\iff$ \ref{casi CI2}: By \cref{rem:casi CI 2}, $a\in \coz L$ iff $\close_L(a)$ is a zero sublocale of $L$, so using the properties of \cite[Section 4]{ferreira_remainders_2018}, we have that 
 \[\close_L(q(a))=\overline{\close_S(q(a))}^L=\overline{\close_L(a) \cap S}^L.\]
Thus, $q(a)\in \coz L$ iff $\overline{\close_L(a) \cap S}^L$ is a zero sublocale of $L$, which shows that $S$ is casi coz-included iff condition (2) holds. \\
 \ref{casi CI1}$\iff$ \ref{casi CI3}: By definition.\\
 \ref{casi CI3}$\iff$\ref{casi CI4}: This is precisely  \cref{casi coz incl iff adj cozeros}.\\
 \ref{casi CI4} $ \implies$ \ref{casi CI5}: If $q_*$ preserves bueno cozero elements, then the restriction $q_*|_{\bcoz S}$ maps to $\coz L$. Furthermore, since $q\dashv q_*$ we have that 
    $$ \coz q (a)= q(a)\leq q(b)\iff a\leq q_*(q(b))=q_*|_{\bcoz S} q(b)$$
for every $a,\,b\in \coz L$. Thus $\coz q$ has a right adjoint. 
\\
\ref{casi CI5}$\implies$ \ref{casi CI6}: Since the right adjoint exists,  by \cref{qcoz q*}, it is the restriction of $q_*$ to $\bcoz S$ (because $\bcoz S = q[\coz L]$).
   By \cref{coz nucleus}, $\upsilon_S$ drops to a nucleus on $\coz L$ iff for each $y\in \coz L$, $\upsilon_S(y)\in \coz L$. Suppose $y \in \coz L$, then $$\upsilon_S(y) = q_*(q(y))=q_*|_{\bcoz S}(q(y))$$(the last equality holds since $q(y)\in \bcoz S= q[\coz L]$). But $\map{q_*|_{\bcoz S}}{\bcoz S}{\coz L}$ and thus $\upsilon_S(y) \in \coz L$.\\
 \ref{casi CI6} $\implies$\ref{casi CI4}: Let $\coz \upsilon_S$ denote the dropped nucleus. Take any bueno cozero $x$ of $S$, then $x= q(y)$ for some $y\in \coz L$. Now $q_*(x) = q_*(q(y))= \upsilon_S(y) = \coz \upsilon_S(y)$ since $y \in \coz L$, and thus by the assumption, this element is in $\coz L$, as desired.\\
 \ref{casi CI7}$ \iff$\ref{casi CI1}: Note that for every $c\in \bcoz S$, since $c\in S$  we have that $c=\tbigvee E_S[c]$. Thus, $\tbigvee E_S[c]\in \coz L$ for every $c\in \bcoz S$ if and only if $\bcoz S\subseteq \coz L$. \\
 \ref{casi CI8}$\iff$\ref{casi CI7}: This holds because $\tbigvee E_S[x]=\tbigvee E_S[q(x)]$ for every $x \in L$.
  \color{black}
 \end{proof}

\begin{remark}\label{rem: casi coz incl}\
\begin{enumerate}[label=\upshape(\arabic*), ref = \theremark\arabic*)]
    \item 
 It is straightforward, using \cref{casi Coz incl charact.}\ref{casi CI2}, to see that every zero sublocale is casi coz-included.\label[remark]{rem: casi coz incl 1} This gives an alternate proof of \cref{when zeros are ci 2}.
  \item Since frame maps preserve cozeros, every frame map ``drops" to a $\sigma$-frame map between the cozeros; a quotient ``drops'' to a quotient iff the frame map is coz-embedded; the associated nucleus ``drops'' iff the frame map is casi coz-included; a quotient ``drops" to a quotient and the nucleus to a nucleus iff the frame map is coz-included. In the latter case, the nucleus is the composition of $\coz q$ and $\coz q_*$.
 \end{enumerate}
 \end{remark}
This weakening of coz-included still suffices to characterize perfect normality (cf. \cref{PN iff all subs coz incl}).
\begin{proposition}\label{prop char pn}
   $L$ is perfectly normal iff every sublocale is casi coz-included iff every closed sublocale is casi coz-included.
\end{proposition}
\begin{proof}
    Assume $L$ is perfectly normal, and let $a\in L$. Any bueno zero sublocale of $\close_L(a)$ is of the form  $Z\cap \close_L(a)$ for some $Z\in\zs(L)$. By assumption $\close_L(a)$ is also a zero sublocale of $L$. Thus, $Z\cap \close_L(a)= \overline{Z\cap \close(a)}^L$ is a zero sublocale of $L$ which proves that $\close_L(a)$ is casi coz-included. 
    On the other hand, assume every closed sublocale is casi coz-included, and let $a\in L$. Trivially, $a\in \bcoz \close_L(a)$. By assumption, $\bcoz \close_L(a)\subseteq \coz L$; thus $a\in \coz L$.
\end{proof}
In fact, as we will see in \cref{char relative zero}\ref{PN} the condition on the sublocales may be further weakened.
\section{Relative zeros}
In this section we introduce a notion that is weaker than both density and casi coz-inclusions. A dense sublocale may be characterized as having the same bottom element as the superlocale, and in many of the proofs above, it is the behavior of the bottom element that is deterministic. Moreover, the bottom element of a zero sublocale is a cozero element (if $S = \close(a)$ for some $a \in \coz L$, then
 $a=0_S \in \coz L$). This leads to the following:
\begin{definition}
     A sublocale $S$ of $L$ is a \emph{relative zero sublocale} if $0_S \in \coz L$. Let $\rzs (L)$ denote all the relative zero sublocales of $L$.
\end{definition}
Equivalently, a frame map is relatively zero if for $\map{h}{L}{M}$, $h_*(0) \in \coz L$. Obviously, any zero sublocale is a relative zero; that is, $\zs (L) \subseteq \rzs (L)$ always.
The naming convention was chosen as it mirrors that of a relatively compact subspace.\footnote{A subset of a topological space is relatively compact if its closure is compact \cite[XI.6]{dugundji_topology_1966}} This is made clear by \cref{relatively zero charact.2}. We need the following well-known facts (see \cite[III.6]{picado_frames_2012} for all the details).

\begin{fact}\label{zero facts} For any $a \in L$,
    \begin{enumerate}[label=\upshape(\arabic*), ref = \thefact\arabic*)]
        \item\label[fact]{zero facts 1}   
         $0_{\close(a)}= a$ and $0_{\open (a)} = a \rightarrow 0 = a^*$.
         \item \label[fact]{zero facts 2}$0_{\close(a)^\circ}=0_{\open (a^*)}= a^* \rightarrow 0 = a^{**}$.
         \item \label[fact]{zero facts 3}$0_{\overline{\open(a)}}=0_{\close (a^*)}= a^*$.
   \end{enumerate}
\end{fact}
\begin{lemma}\label{open and closed as relative zeros}\
   \begin{enumerate}[label=\upshape(\arabic*), ref = \thelemma\arabic*)]
        \item A closed sublocale $\close(a)$ is a relative zero sublocale iff $a \in \coz L$ iff it is a zero sublocale;
         that is, $$\close(a) \in \rzs(L) \iff a \in \coz L \iff \close(a) \in \zs(L).$$
        \item An open sublocale $\open(a)$ is a relative zero  sublocale iff $a^* \in \coz L$ iff the closure of $\open (a)$ is a zero sublocale\label[lemma]{open and closed as relative zeros 2}; that is,
        $$\open(a) \in \rzs(L) \iff a^* \in \coz L \iff \overline{\open(a)}\in \zs(L). $$
         \end{enumerate}
\end{lemma}
\begin{proof} These are straightforward consequences of \cref{zero facts}.
    \end{proof}
\begin{definition}An element $a \in L$ is a \emph{relative zero} if $\open(a)$ is a relative zero sublocale; that is, $a^* \in \coz L$. Let  $\rz L$ denote the set of all relative zero elements of $L$. 
\end{definition} 
\begin{remark}\
    \begin{enumerate}
        \item Note that for any sublocale $S$, $0_S=0_{\overline{S}}$, and thus $S$ is a relative zero sublocale if and only if $\overline{S}$ is a zero sublocale.\footnote{The taking of closures of zeros has been considered in point-set topology - see \cite[6E]{gillman_rings_1960}}
        \item Since the closure of a finite join of sublocales is equal to the join of their closures, the finite join of relative zero sublocales is relative zero. That is, $\rzs (L)$ is a sub-join semi-lattice of $\Ss(L)$.
    \end{enumerate}
\end{remark}
For the sake of comparisons to \cref{Coz incl charact.,casi Coz incl charact.} we include the following equivalences. The proofs are straightforward.      
        \begin{theorem}\label{relatively zero charact.} For any sublocale $S$ of $L$, the following are equivalent:
 \begin{enumerate}[label=\upshape(\arabic*), ref = \thetheorem(\arabic*)]
     \item $S$ is a relative zero sublocale.
     \item $\overline{S}$ is a zero sublocale of $L$\label[theorem]{relatively zero charact.2}.
     \item The quotient $\map{q}{L}{S}$ is relatively zero.
     \item For the right adjoint, $q_* (0_S) \in \coz L$.
      
      \item For the associated nucleus, $\upsilon_S(0_L)\in \coz L$.
      \item For the associated congruence, $\bigvee E_S[0_L] \in \coz L$.
 \end{enumerate}     
\end{theorem}

       The next proposition establishes the relation between relative zero sublocales and other special sublocales of interest in this paper.
\begin{proposition}\label{relatively zero implications}
    For a sublocale $S \subseteq L$,
    \begin{enumerate}[label=\upshape(\arabic*),ref = \theproposition(\arabic*)]
    \item Dense implies relative zero, but not conversely.
        \item Casi coz-included implies relative zero, but not conversely.
        \item Coz-embedded and relative zero are independent notions.\label[proposition]{relatively zero implications 3}
    \end{enumerate}
\end{proposition}
\begin{proof} Both forward implications of (1) and (2) are straightforward.  For (1): any proper zero sublocale is relative zero but not dense.
    For (2): any dense sublocale that is not casi coz-included; for example the Booleanization of $\beta \mathfrak L \Rn$ (see \cref{example of casi not coz included}).\\
    (3) Consider a normal, not perfectly normal frame (for example, $\beta \mathfrak L \Rn$), then there exists a closed, non zero sublocale, say $\close(a)$ for $a \notin \coz L$.
    Since $0_{\close(a)} =a$, this sublocale is not relatively zero, but it is coz-embedded (\cref{char. of normality}). On the other hand, by \cref{example of casi not coz included} the frame $\mathcal{O}E$ is not coz-embedded in $\mathcal{O}X$, but it is relative zero (since it is dense).
\end{proof}
\begin{remark}
    The first example in \cref{relatively zero implications 3} shows that both $C$ and $C^*$-embedded sublocales need not be relative zeros.
\end{remark}

Coz-inclusions may be characterized via relative zero sublocales:
\begin{proposition}
    Let $S$ be a sublocale of $L$, then the following are equivalent:
    \begin{enumerate}[label=\upshape(\arabic*)]
        \item\label{ci} $S$ is coz-included. 
        \item\label{rz.rz} Every  relative zero sublocale of $S$ is a relative zero sublocale of $L$.
        \item\label{z.rz} Every zero sublocale of $S$ is a relative zero sublocale of $L$.
    \end{enumerate}
    That is, 
    $$\coz S \subseteq \coz L \iff \rzs (S) \subseteq \rzs (L) \iff \zs (S) \subseteq \rzs (L).$$
\end{proposition}

\begin{proof}
    \ref{ci}$\implies$\ref{rz.rz}: Let $T$ be a relative zero sublocale of $S$. That is, $0_T\in \coz S$. By assumption $\coz S\subseteq \coz L$; thus, $0_T\in \coz L$ which means $T$ is relative zero in $L$. \\
    Since every zero sublocale of $S$ is relative zero in $S$, \ref{rz.rz} trivially implies \ref{z.rz}.\\
    \ref{z.rz}$\implies$\ref{ci}: Using \cref{{relatively zero charact.2}}, this is precisely implication \ref{CI:part2}$\implies$\ref{CI:part1} of Theorem~\ref{Coz incl charact.}.
\end{proof}
The obvious weakening of the above yields the following characterization of casi coz-inclusions.
\begin{proposition}
    Let $S$ be a sublocale of $L$, then the following are equivalent:
    \begin{enumerate}[label=\upshape(\arabic*)]
        \item \label{cci} $S$ is casi coz-included.
        \item \label{bz.rz} Every bueno zero sublocale of $S$ is a relative zero of $L$.
    \end{enumerate}
    That is, $$\bcoz S \subseteq \coz L \iff \bzs (S)\subseteq \rzs (L).$$
\end{proposition}

The characterization of when open sublocales are relatively zero, that is, $a^* \in \coz L$ (see \cref{open and closed as relative zeros}) reminds one of certain topological properties relating cozeros and pseudocomplements. 
The notions of Oz and weakly Oz frames are already known concepts (both in spaces and frames). The third definition below is a new class obviously relevant from our perspective. We include definitions/characterizations of the former two concepts to make the comparisons more transparent.  Note that, in general, $(\coz L)^* \not= (\coz L)^{**}$ even though $L^* = L^{**}$. 
\begin{definition}
    A completely regular frame $L$ is 
    
    \begin{enumerate}[label=\upshape(\arabic*), ref = \thedefinition(\arabic*)]
        \item \cite[Prop.~2.2]{banaschewski_oz_2009} \emph{$Oz$} if every regular element is a cozero;\\ that is, $L^{**} \subseteq \coz L$;
        \item \cite[Defn.~5.1]{banaschewski_oz_2009} \emph{weak Oz} if $a^*\in \coz L$ for each $a\in \coz L$;\\ that is, $(\coz L)^* \subseteq \coz L$;
        \item \emph{casi Oz} if if $a^{**}\in \coz L$ for each $a\in \coz L$; that is, $(\coz L)^{**} \subseteq \coz L$.
\end{enumerate}
\end{definition}
In \cref{Coole vs Oz 2} we will give a condition that renders all three equivalent.
\begin{remark}\
    
\begin{enumerate}[label=\upshape(\arabic*), ref = \theremark(\arabic*)]
\item $L$ is Oz iff every $a \in L$ is coz-embedded \cite[Defn.~2.1]{banaschewski_oz_2009}.\footnote{This is taken as the definition in \cite{banaschewski_oz_2009}, and then the definition given above is a characterization.}
    \item Weak Oz is, as one would expect, weaker than Oz, and strictly so: for any Boolean $\sigma$-algebra, the frame of its ideals is weak Oz, but it is Oz iff the Boolean algebra is complete. See \cite[Ex.~5.1]{banaschewski_oz_2009}.
 \item Casi Oz is weaker than weak Oz and also strictly so: the one point compactification of an uncountable set with the discrete topology (see \cref{1-pt comp} ) is casi Oz but not weak Oz. Moreover, the use of ``casi" as an adjective of Oz fits with the philosophy that this should be equivalent when dealing with coz-embeddings (as will be seen in \cref{Coole vs Oz 2}).
 \item Clearly both perfect normality and extreme disconnectivity imply Oz. Recall that a  frame is \emph{basically disconnected}  if $CL = (\Coz L)^{**}$ \cite[8.4.2]{ball_c-_2002} and hence basically disconnected implies casi Oz; in fact it implies weak Oz.\footnote{It seems that this rather obvious fact was overlooked in \cite{banaschewski_oz_2009}.} The latter follows since for basically disconnected frames, $(\Coz L)^{**}=(\Coz L)^*$. 
 \end{enumerate}
 \end{remark}
Oz frames are characterized as having every open sublocale  coz-embedded. A slightly stronger condition would be that every open sublocale is coz-included (perhaps these should be ``perfect Oz'' frames); a condition that is formally weaker than perfect normality. We do not yet have a differentiating example to prove that the condition is strictly weaker, nor have we been able to prove that the two notions are equivalent.

We now consider characterizations of the weaker variants of perfect normality (including perfect normality) introduced above in terms of relative zeros. We recall that $a^{***}= a^*$, and $0_{\open (a)} = a^*$ for any $a \in L$ ( \cref{zero facts}(1)).
\begin{proposition}\label{char relative zero}\
\begin{enumerate}[label=\upshape(\arabic*)]
     \item\label{casi Oz} $L$ is casi Oz iff the interior of every zero sublocale is relatively zero iff every pseudocomplement of a cozero element is relatively zero. 
That is, $$(\coz L)^{**} \subseteq \coz L \iff \zs^\circ(L)\subseteq \rzs(L)\iff (\coz L)^{*} \subseteq \rz L$$
where $\zs^\circ(L)=\{ Z^\circ \mid Z\in \zs(L)\}$.
    \item \label{weak Oz} $L$ is weak Oz iff every cozero sublocale is relatively zero iff every cozero element is relatively zero. That is, $$(\coz L)^{*} \subseteq \coz L \iff \cozs(L)\subseteq\rzs(L)\iff \coz L \subseteq \rz L.$$ 
    \item \label{Oz} $L$ is Oz iff every open sublocale is relatively zero iff every element of $L$ is relatively zero. That is, 
    $$ L^{**} \subseteq \coz L \iff \os(L)\subseteq\rzs(L)\iff L = \rz L.$$ 
    \item \label{PN} $L$ is perfectly normal iff every closed sublocale is relatively zero iff every element of $L$ is cozero. That is, 
    $$L = \coz L \iff \cs(L) = \rzs (L).$$ 
\end{enumerate}
\end{proposition}

\begin{proof}
\ref{casi Oz} By definition, $L$ is casi Oz iff $(\coz L)^{**} \subseteq \coz L$. The following equivalences give the required results:
\begin{align*}
  a^* \in \rz (L) \mbox{ for } a \in \coz L     &\iff  a^{**}\in \coz L \mbox{ for } a \in \coz L \quad\mbox{ ($\Leftarrow$ uses $a^{***}=a^*$)} \\ &\iff O_{\close(a)^\circ} \in \coz L \mbox{ for } a \in \coz L \quad\mbox{(by \cref{zero facts}(2)).}\\
\end{align*}
\ref{weak Oz} For any $c \in \coz L$, $c^* \in \coz L$ iff  $0_{\open (c)} \in \coz L$.\\
\ref{Oz} Let $a \in L$. Since $a^* \in L^{**}$, if $L$ is Oz, then $a^* \in \coz L$, so $a$ is a relative zero. Conversely, since $a^*$ is a relative zero, $a^{**}\in \coz L$, showing that $L$ is Oz.\\ 
\ref{PN} The forward implication is clear. For the reverse: suppose $a \in L$, then  $0_{\close (a)} \in \coz L$, but $a = 0_{\close (a)}$.

\end{proof}

Recall that a frame is coole if the Booleanization map is coz-onto (\cref{defn: coole}). In \cite[Prop.~5.2]{banaschewski_oz_2009} a condition is introduced to force the following equivalence: a weak Oz frame $L$ is Oz iff every regular element of $L$ is $c^*$ for some $c\in \Coz L$, that is $L^{**} \subseteq (\Coz L)^*$. Various names have been used for this condition: Greunhage \cite{gruenhage_products_2006} chose the term ``cozero approximated'' to describe the same property for topological spaces; Hager and Martinez \cite{hager_fraction-dense_1993} used ``fraction dense''.
This condition is equivalent to the frame being coole.
\begin{proposition}
$L^{**} \subseteq (\coz L)^*$ iff $L$ is coole.
\end{proposition}
\begin{proof} Note that $\coz ( L^{**}) =  L^{**}$ since all elements are complemented, and hence cozero. Suppose $L^{**} \subseteq (\coz L)^*$. Take any $x\in L^{**}$, then $x^*\in L^{**}$, so by assumption $x^*=c^*$ for some cozero $c$, and thus $x=x^{**}=c^{**}$, which shows the required coz-ontoness. 

Conversely, suppose $L \to L^{**}$ is coz-onto. Take any $x \in L^{**}$ then $x^* \in L^{**}$, so by assumption $x^*=c^{**}$ for some $c \in \coz L$, and thus $x =x^{**}=c^{***}=c^*$.

\end{proof}
\begin{samepage}
\begin{corollary}\label{Coole vs Oz}\
\begin{enumerate}[label=\upshape(\arabic*), ref = \thecorollary(\arabic*)]
    \item All Oz frames are coole, but not conversely.

\item For coole frames, all three notions coincide. That is, a coole frame is Oz iff it is weak Oz iff it is casi Oz.\label[corollary]{Coole vs Oz 2}
\end{enumerate}
\end{corollary}
\end{samepage}
\begin{proof}(1) For any $L$, $\coz ( L^{**}) =  L^{**}$, and thus the result follows since $L$ being Oz means $L^{**}\subseteq \coz L$. 
For the converse, see \cref{beta R is not Oz} and \cref{betaR is coole}: $\beta \mathfrak L \Rn$ is coole but not Oz.\\
(2) The first equivalence is the result mentioned above \cite[Prop.~5.2]{banaschewski_oz_2009}. Suppose $L$ is coole and casi Oz.
   Take $x\in L$ then, by cooleness, $x^{**}=c^{**}$ for some $c\in \coz L$, but casi Oz gives that $c^{**} \in \coz L$, and thus so is $x^{**}$. 
\end{proof}    
 
\begin{remark}\
    \begin{enumerate}[label=\upshape(\arabic*), ref = \theremark(\arabic*)]
     \item From a topological viewpoint, one can deduce that Oz spaces are coole, that is, the Booleanization of the topology is coz-embedded. This was not mentioned in \cite{gruenhage_products_2006}, although perfectly normal and ED spaces (both subsets of Oz) were noted to be coole.
    \item Whenever the Booleanization is Lindel\"of, the frame is coole: every quotient with Lindel\"of codomain is a coz-embedding (see \cite[prop.~3.2]{dube_coz-onto_2006}). Such frames are precisely the $q$-Lindel\"of ones (see \cite[defn.~3.1]{dube_when_2017} (since $q$-Lindel\"of is equivalent to the Booleanization being Lindel\"of \cite[prop.~3.6]{dube_when_2017}), and thus $q$-Lindel\"of implies coole. 
\end{enumerate}
\end{remark}

\begin{proposition} For a completelely regular frame $L$,
 $L$ is casi Oz iff the Booleanization is casi coz-included.
\end{proposition}
\begin{proof}By definition, the Booleanization is casi coz-included iff $(\coz L)^{**} \subseteq \coz L$.
    \end{proof}

Finally we have sufficient machinary to examine, and mend, the second of the errors  in \cite{ball_lindelof_2017} which we refered to above:
 \cite[Prop.~6.2]{ball_lindelof_2017} is incorrect as it is stated (as is the proof since coz-ontoness is NOT inherited by superlocales). However with the slightly stronger notion of coz-included the following is true: \cref{Coole vs Oz 2} and the discussion above provide the correct proof. 
\begin{proposition}\label{char of Oz}
 The following are equivalent:
 \begin{enumerate}[label=\upshape(\arabic*), ref = \theproposition(\arabic*)]
    \item  $L$ is an Oz frame.
 \item Every pseudocomplement is cozero.
\item  Every open sublocale is coz-embedded.
\item  Every dense sublocale is coz-embedded.
 \item  The Booleanization $L^{**}$ is coz-included.
 \item   $L$ is  coole and casi Oz.
\end{enumerate}
\end{proposition}

\begin{corollary}\label{beta R is not Oz}
    
The \v{C}ech-Stone compactification of the locale of the reals, $\beta (\mathfrak{L} \Rn)$, is not Oz.
    
\end{corollary}
\begin{proof}
    We saw in \cref{ex. coole not Oz} that $(\mathfrak{L} \Rn)^{**}$ is not coz-included in $\beta (\mathfrak{L} \Rn)$.
   
\end{proof}
Since Oz is a conservative notion \cite{banaschewski_oz_2009} and compact regular frames are spatial, this is an alternate, and somewhat less complicated, proof of the result by Terada \cite{terada_spaces_1979}  showing that that $\beta \Rn$ is not Oz.

\section{In summation}
 We recall the following definitions: a completely regular frame is
 
 \begin{itemize}
 \item a \emph{P-frame} if every cozero element is complemented \cite[Defn.~8.4.6]{ball_c-_2002};
       \item \emph{almost Boolean } if it is an extremally disconnected P-frame \cite{banaschewski_nonmeasurable_2015}.
        \end{itemize}
\cref{Booleanization table} characterizes frames having their Booleanization  embedded in a designated way:

\begin{table}[h]
    \centering
    \begin{tabular}{|c|c|}
        \hline
       \textbf{Booleanization} & \textbf{Frame}  \\ \hline\hline
       coz-embedded & coole \\ \hline
       casi coz-incl. &casi-Oz\\ \hline
       coz-included & Oz \\ \hline
        $C^*$-embedded &ED\\ \hline
        $C$-embedded &almost Boolean \\ \hline
        isomorphism &Boolean \\ \hline
    \end{tabular}
    \caption{}
    \label{Booleanization table}
\end{table}

Various topological properties have been characterized in terms of designated sublocales (or subspaces) being embedded in a specific way - perhaps the most well known is that for closed sets. Our aim is to give a summary of known results incorporating the  new additions that were developed above, together with the following:

 \begin{proposition} For a completely regular frame $L$,
\begin{enumerate}[label=\upshape(\arabic*), ref = \theproposition(\arabic*)]
    \item  $L$ is almost Boolean iff all open sublocales are C-embedded.
    \item $L$ is a weakly perfectly normal $P$-frame iff all sublocales are C-embedded.
\end{enumerate}
\end{proposition}
\begin{proof}
    (1) Let $L$ be an almost Boolean frame. Since $L$ is extremally disconnected, every open sublocale is $C^*$-embedded (\cite[Proposition 8.4.1]{ball_c-_2002}). In particular, every open sublocale is coz-embedded. Recall from \cite[Proposition 4.9]{dube_coz-onto_2006} that in a $P$-frame every coz-embedded sublocale is $C$-embedded. Thus, every open sublocale is $C$-embedded.
    Conversely, if all opens are $C$-embedded, in particular all cozero sublocales are $C$-embedded and all open sublocales are $C^*$-embedded. Thus, by Propostion 8.4.7 and Proposition 8.4.1 of  \cite{ball_c-_2002}, $L$ is a $P$-frame and it is extremally disconnected.\\
    (2) Let $L$ be a weakly perfectly normal $P$-frame. Recall that a frame is weakly perfectly normal if all sublocales are coz-embedded (see \cref{sec coz inclusions}). Suppose $S$ is a sublocale of $L$, then it is coz-embedded. Moreover, in a $P$-frame every coz-embedded sublocale is $C$-embedded (\cite[Proposition 4.9]{dube_coz-onto_2006}), so $S$ is $C$-embedded, as required.
    Conversely, if all sublocales are $C$-embedded, then all cozeros are $C$-embedded. Thus, $L$ is a $P$-frame. Furthermore, all sublocales are coz-embedded ($C$-embedded implies coz-embedded) so $L$ is also weakly perfectly normal. 
\end{proof}
\begin{corollary}
    A weakly perfectly normal $P$-frame is extremally disconnected, hence almost Boolean. 
\end{corollary}
\begin{remark}
   We have seen that Oz is stronger than coole, but just as Oz + almost P (iff ED +P) iff almost boolean, it is straightfoward to show that coole + almost P iff almost Boolean.
\end{remark}

We list some definitions that we have not yet considered in this paper:\\
 A completely regular frame $L$ is 
\begin{itemize}
    \item \emph{completely normal} if for every $a,b \in L$, there exists disjoint $x, y \in L$ with $b \leq a \join y$ and $a \leq b \join x$ (see \cite[p.~6]{ferreira_completely_2009});\footnote{Equivalently, every pair of separated sublocales of $L$ is separated by open sublocales.}
    
     \item an \emph{F-frame } if every cozero is $C^*$-embedded, and a \emph{quasi F-frame } if every dense cozero is $C^*$-embedded \cite[Defn.~8.4.9]{ball_c-_2002}; 
      \item an \emph{almost P-frame}\footnote{Also known as $\omega_1$-hollow frames - see \cite{ball_lambda-hollow_2024}.} if $a = a^{**}$ for all $a \in \coz L$  \cite[Defn.~8.4.6]{ball_c-_2002}.
 \end{itemize}

We also highlight some results and facts that are required to understand \cref{Sublocale table}:\\
 For a completely regular frame $L$ 
\begin{itemize}
    \item Every cozero sublocale is coz-embedded (see \cite[3.2.11]{ball_c-_2002}).
   \item A zero sublocale is coz-embedded iff it is  $C^*$-embedded iff it is C-embedded (see \cite[Cor.~3.4.6]{avilez_point-free_2022}) iff it is coz-included (\cref{when zeros are ci 2}). This explains the $\star$ in \cref{Sublocale table}.
   \item Extremally disconnected and completely normal is equivalent to all sublocales are $C^*$-embedded \cite[Prop.~6.7]{ferreira_completely_2009}.
   \item A frame is a P-frame iff all cozero sublocales are $C$-embedded (see \cite[Prop.~8.4.7]{ball_c-_2002} and \cite[4.3.19]{avilez_point-free_2022}).
   \item Since casi coz-included implies relative zero, if every open is casi coz-included, then the frame is Oz (see \cref{char relative zero} \ref{Oz}). Hence, all open sublocales are casi coz-included iff  they are coz-included; that is, the frame is ``perfectly Oz''.
\end{itemize}

 \cref{Sublocale table} characterizes those frames where all designated sublocales (dense, cozero zero, open, closed or all) are embedded in a specified way (coz-embedded, $C^*$-embedded, $C$-embedded, relative zero, casi coz-included or coz-included). We use the following abbreviations:
\begin{itemize}
    \item[] ED: extremally disconnected
    \item[] PN: perfectly normal
     \item[] WPN: weakly perfectly normal
     \item[] almB: almost Boolean
     \item[] POz: perfectly Oz
     \item[] $\star, \dagger, \ast , \circ$ denote equivalences of the collection of sublocales, none of which have been successfully characterized by a known notion.
\end{itemize}

\begin{table}[h]
    \footnotesize
    \centering
    \begin{tabular}{|c|c|c|c|c|c|c|c|c|}
        \hline
       \bf{Sublocales} &\bf{dense Cozero}& \bf{Cozero} & \bf{Zero} & \bf{Open} & \bf{Closed} & \bf{Dense} & \bf{All} \\\hline\hline
       
        coz-emb. & all frames& all frames &\color{blue}{$\star$} & Oz & Normal & Oz  & WPN \\ \hline
        $C^*$-emb.&quasi F &F-frame &\color{blue}{$\star$} &ED &Normal &ED & {ED+CN} \\ \hline
        $C$-emb. &almost P&P-frame & \color{blue}{$\star$}&{almB}&Normal &almB &{P + WPN} \\ 
        \hline\hline
        
        rel. zero &all frames &weak Oz& all frames &{Oz}&{PN}&all frames&PN\\ \hline
        casi coz-incl.&\color{blue}{$\ast$}&\color{blue}{$\dagger$} &all frames  & {POz}&{PN} &\color{blue}{$\circ$} & {PN}\\ \hline
        coz-incl.&\color{blue}{$\ast$} &\color{blue}{$\dagger$} &\color{blue}{$\star$} &{POz} & PN&\color{blue}{$\circ$+Oz}&PN \\ \hline
        
    \end{tabular} 
    \caption{}
    \label{Sublocale table}
\end{table}

\bibliographystyle{alpha}

\begin{thebibliography}{BMPWW17}
 
    \bibitem[Avi22]{avilez_point-free_2022}
A.~B. Avilez.
\newblock {\em A {Point}-{Free} {Study} of z-{Embeddings}, {More} {General}
  {Classes} of {Localic} {Maps}, and {Uniform} {Continuity}}.
\newblock PhD thesis, Universidade de Coimbra, December 2022.

\bibitem[Ban97]{banaschewski_real_1997}
B.~Banaschewski.
\newblock {\em The {Real} {Numbers} in {Pointfree} {Topology}}.
\newblock Universidad de Coimbra, 1997.

\bibitem[Ban15]{banaschewski_nonmeasurable_2015}
B.~Banaschewski.
\newblock Nonmeasurable cardinals and pointfree topology.
\newblock {\em Math. Slovaca}, 65(2):289--300, 2015.

\bibitem[BDGWW09]{banaschewski_oz_2009}
B.~Banaschewski, T.~Dube, C.~Gilmour, and J.~Walters-Wayland.
\newblock Oz in pointfree topology.
\newblock {\em Quaest. Math.}, 32(2):215--227, 2009.

\bibitem[BG89]{banaschewski_stonecech_1989}
B.~Banaschewski and C.~Gilmour.
\newblock Stone-{\v{C}ech} {compactification} and {dimension} {theory} for
  {regular} σ-{frames}.
\newblock {\em J. Lond. Math. Soc.}, s2-39(1):1--8, 1989.

\bibitem[BG01]{banaschewski_cozero_2001}
B.~Banaschewski and C.~Gilmour.
\newblock Cozero bases of frames.
\newblock {\em J. Pure Appl. Algebra}, 157(1):1--22, 2001.

\bibitem[BH74]{blair_extensions_1974}
R.~L. Blair and A.~W. Hager.
\newblock Extensions of zero-sets and of real-valued functions.
\newblock {\em Math. Z.}, 136(1):41--52, 1974.

\bibitem[BHWW24]{ball_lambda-hollow_2024}
R.~N. Ball, A.~W. Hager, and J.~Walters-Wayland.
\newblock From $\lambda$-hollow frames to $\lambda$-repletions in {W} {I}: {The}
  scaffolding of a completely regular frame.
\newblock {\em Topology Appl.}, 342:108774, 2024.

\bibitem[Bla76]{blair_spaces_1976}
R.~L. Blair.
\newblock Spaces in which special sets are z-embedded.
\newblock {\em Can. J. Math.}, 28:673--690, 1976.

\bibitem[BMPWW17]{ball_lindelof_2017}
R.~N. Ball, M.~A. Moshier, A.~Pultr, and J.~L. Walters-Wayland.
\newblock Lindelöf tightness and the {Dedekind}-{MacNeille} completion of a
  regular σ-frame.
\newblock {\em Quaest. Math.}, 40(3):347--362, 2017.

\bibitem[BWW02]{ball_c-_2002}
R.~N. Ball and J.~Walters-Wayland.
\newblock C- and {C} $^{{*}}$-quotients in pointfree topology.
\newblock {\em Dissertationes Math.}, 412:1--62, 2002.

\bibitem[BWW06]{ball_well-embedding_2006}
R.~N. Ball and J.~L. Walters-Wayland.
\newblock Well-embedding and {G}$_\delta$-{density} in a {pointfree} {setting}.
\newblock {\em Appl. Categ. Structures}, 14(4):351--355, 2006.

\bibitem[Dub17]{dube_when_2017}
T.~Dube.
\newblock When {Boole} commutes with {Hewitt} and {Lindelöf}.
\newblock {\em Appl. Categ. Structures}, 25(6):1097--1111, 2017.

\bibitem[Dug66]{dugundji_topology_1966}
J.~Dugundji.
\newblock {\em Topology}.
\newblock Series in {Advanced} {Math}. Allyn and Bacon, Boston, 1966.

\bibitem[DWW06]{dube_coz-onto_2006}
T.~Dube and J.~Walters-Wayland.
\newblock Coz-onto frame maps and some applications.
\newblock {\em Appl. Categ. Struct.}, 15(1):119--133, 2006.

\bibitem[FGGP09]{ferreira_completely_2009}
M.~J. Ferreira, J.~Gutiérrez~García, and J.~Picado.
\newblock Completely normal frames and real-valued functions.
\newblock {\em Topology Appl.}, 156(18):2932--2941, 2009.

\bibitem[FPP18]{ferreira_remainders_2018}
M.~J. Ferreira, J,~ Picado, and S.~M.~ Pinto.
\newblock Remainders in pointfree topology.
\newblock {\em Topology Appl.}, 245:21--45, August 2018.

\bibitem[GGPP14]{gutierrez_garcia_notes_2014}
J.~Gutiérrez~García, J.~Picado, and A.~Pultr.
\newblock Notes on point-free real functions and sublocales.
\newblock In {\em Categorical {Methods} in {Algebra} and {Topology}: {Special}
  {Volume} in {Honour} of {Manuela} {Sobral}}, volume~46, pages 167--200. DMUC
  - Textos de Matemática, 2014.

\bibitem[GJ60]{gillman_rings_1960}
L.~Gillman and M.~Jerison.
\newblock {\em Rings of continuous functions}.
\newblock D. Van Nostrand Company, Inc., Princeton-Toronto-London-New York,
  1960.

\bibitem[Gru06]{gruenhage_products_2006}
G.~Gruenhage.
\newblock Products of cozero complemented spaces.
\newblock {\em Houston J. Math.}, 32(3):757--773, 2006.

\bibitem[HM93]{hager_fraction-dense_1993}
A.~W. Hager and J.~Martinez.
\newblock Fraction-{dense} {algebras} and {spaces}.
\newblock {\em Canad. J. Math.}, 45(5):977--996, 1993.

\bibitem[PP12]{picado_frames_2012}
J.~Picado and A.~Pultr.
\newblock {\em Frames and {Locales}: {Topology} without points}.
\newblock Frontiers in {Mathematics}. Springer, Basel, 2012.

\bibitem[Ter79]{terada_spaces_1979}
T.~Terada.
\newblock On spaces whose {Stone}-Čech compactification is {Oz}.
\newblock {\em Pacific J. Math.}, 85(1):231--237, 1979.

\end{thebibliography}

\end{document}